\documentclass[10pt]{amsart}
\usepackage{latexsym}
\usepackage{amsmath}
\usepackage{amssymb}
\usepackage[all]{xy}

\numberwithin{equation}{section}
\newtheorem{thm}{Theorem }[section]

\newtheorem{prop}[thm]{Proposition}

\newtheorem{defi}[thm]{Definition}

\newtheorem{ex}[thm]{Example}
\newtheorem{qu}[thm]{Question}
\newtheorem{rem}[thm]{Remark}
\newtheorem{prob}[thm]{Problem}

\newcommand{\C}{{\mathbb C}}
\newcommand{\Q}{{\mathbb Q}}

\begin{document}
\title{Rational toral rank of a map
}

\author{Toshihiro YAMAGUCHI }
\address{Faculty of Education, Kochi University, 2-5-1,Kochi,780-8520, JAPAN}
\email{tyamag@kochi-u.ac.jp}
\maketitle

\begin{abstract}
Let $X$ and  $Y$ be   simply connected   CW complexes 
with finite rational cohomologies.
The rational toral rank $r_0(X)$ of a space $X$
is the largest integer $r$ such that the torus
 $T^r$  can act continuously
 on a CW-complex  in  the rational homotopy type of $X$
 with all its isotropy subgroups finite \cite{H}.
As a rational  homotopical condition to be a toral map  preserving   almost free toral actions
for a map $f:X\to Y$,
we define the rational toral rank $r_0(f)$ of $f$, which is a natural invariant
with $r_0(id_X)=r_0(X)$ for the identity map  $id_X$ of $X$.
We will see some properties of it by Sullivan models, which is a free commutative differential graded algebra over $\Q$
\cite{FHT}.
\end{abstract}

\section{Introduction}

We assume that  spaces $X$ and  $Y$  are  simply connected   CW complexes 
with finite rational cohomologies.
Let $T^r$ be  an $r$-torus
 $S^1 \times\dots\times S^1$($r$-factors)
and let $r_0(X)$ be the {\it rational  toral rank}, which is 
 the largest integer $r$ such that a
 $T^r$  can act continuously
 on a CW-complex  in  the rational homotopy type of $X$
 with all its isotropy subgroups finite \cite{H}.
Such an action is called  {\it almost free}. 
Our motivation is in the following 
problem for an equivariant
property of  a map $f:X\to Y$.
 
\begin{prob}\label{1.1}   For an  almost free $T^r$-action   $\mu$ on $X$,
 when  can one put    an  almost free  $T^r$-action    on $Y$
so that $f$ becomes $T^r$-equivariant  ? 
Conversely,  given  an   almost free $T^r$-action $\tau$  on $Y$,
when does $X$ admit an almost free $T^r$-action  making $f$ an  $T^r$-equivariant
map ?
  $$ 
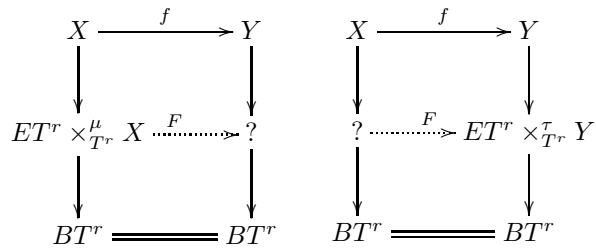
\xymatrix{ X\ar[d]\ar[r]^f&Y\ar[d]\\
  ET^r\times_{T^r}^{\mu}X\ar[d]\ar@{.>}[r]^{\ \ F}&?\ar[d]\\
  BT^r\ar@{=}[r]&BT^r
 }
\ \ \ \ 
\xymatrix{ X\ar[d]\ar[r]^f&Y\ar[d]\\
?  \ar[d]\ar@{.>}[r]^{F\ \ \ }&ET^r\times_{T^r}^{\tau}Y\ar[d]\\
 BT^r\ar@{=}[r]&BT^r
  }$$
\end{prob}
 
Here $X\to  ET^r\times_{T^r}^{\mu}X\to BT^r$ means the Borel fibration
of  a $T^r$-action $\mu$ on $X$.
The integer $r$ of Problem \ref{1.1} is   bounded from above  by the following
numerical invariant,
obtained from a  diagram which  is a rational homotopy version of a $T^r$-equivariant map
for almost free $T^r$-actions.
In this paper, we propose
 
\begin{defi}\label{11} {\rm
For a map $f:X\to Y$, we say  that  the {\it rational toral rank of $f$},
denoted as  $r_0(f)$,    
 is $r$ when it is  the largest  integer  such that there is  a map $F$ between 
fibrations over $BT^r_{\Q}$: 
$$ \ \xymatrix{X_{\Q}\ar[d]_i\ar[r]^{f_{\Q}}& Y_{\Q}\ar[d]^i\\
E_1\ar[d]_p\ar@{.>}[r]^{F} &E_2\ar[d]^p\\
BT^r_{\Q}\ar@{=}[r]&BT^r_{\Q}
}\ \ \ \ \ (*)$$
with $\dim H^*(E_i;\Q)<\infty$ for $i=1,2$.
}\end{defi}
Here $X_{\Q}$ and $f_{\Q}$ are  the rationalizations  \cite{HMR} of a 
simply connected CW complex  $X$
of finite type 
and a map $f$, respectively.
Let 
the  Sullivan minimal model of $X$ be $M(X)=(\Lambda {V},d)$. 
  It is a free $\Q$-commutative differential graded algebra (DGA)
 with a $\Q$-graded vector space $V=\bigoplus_{i\geq 2}V^i$
 where $\dim V^i<\infty$ and a decomposable differential; i.e., $d(V^i) \subset (\Lambda^+{V} \cdot \Lambda^+{V})^{i+1}$ and $d \circ d=0$.
 Here  $\Lambda^+{V}$ is 
 the ideal of $\Lambda{V}$ generated by elements of positive degree. 
Denote the degree of a homogeneous element $x$ of a graded algebra as $|{x}|$.
Then  $xy=(-1)^{|{x}||{y}|}yx$ and $d(xy)=d(x)y+(-1)^{|{x}|}xd(y)$. 
Note that  $M(X)$ determines the rational homotopy type of $X$.
In particular,  $H^*(\Lambda {V},d)\cong H^*(X;\Q )$
and $V^i\cong Hom(\pi_i(X),\Q)$.
Refer to  \cite{FHT} for details.
If an $r$-torus $T^r$ acts on a simply connected space $X$
by $\mu :T^r\times X\to X$, there is the Borel fibration
$$
X \to ET^r \times_{T^r}^{\mu} X \to BT^r,
$$
where 
$ ET^r \times_{T^r}^{\mu} X $  is the orbit space of the  action
$g(e,x)=(e\cdot g^{-1},g\cdot x)$  
on the product $ ET^r \times  X $ for any $e\in ET^r$, 
$x\in X$ and $g\in T^r$.
Note that $ET^r \times_{T^r}^{\mu} X$ 
is rational  homotopy equivalent  to the $T^r$-orbit space of $X$
when $\mu$ is an  almost free toral action \cite{FOT}.
The above Borel fibration  is rationally given by the relative model (Koszul-Sullivan (KS) model)
$$
(\Q[t_1,\dots,t_r],0)
 \to (\Q[t_1,\dots,t_r] \otimes \Lambda {V},D)
 \to (\Lambda {V},d)\ \ \ \ \ (**)$$
where    with $|{t_i}|=2$ for $i=1,\dots,r$, $Dt_i=0$ and
$Dv \equiv dv$ modulo the ideal $(t_1,\dots,t_r)$ for $v\in V$.
The following 
 criterion of Halperin
is  used in this paper.
\begin{prop}\cite[Proposition 4.2]{H}\label{H}
Suppose that $X$ is a simply connected CW-complex  with 
$\dim H^*(X;\Q)<\infty$.
Put $M(X)=(\Lambda V,d)$.
Then  $r_0(X) \ge r$ if and only if there is a relative  model  $(**)$
 satisfying $\dim H^*(\Q[t_1,\dots,t_r] \otimes \Lambda {V},D)<\infty$.
 Moreover,
 if  $r_0(X) \ge r$,
 then $T^r$ acts freely on a finite complex $X'$
 in the  rational homotopy type of  $X$
 and $M(ET^r\times_{T^r}X')\cong 
(\Q[t_1,\dots,t_r] \otimes \Lambda {V},D)$.
\end{prop}

The diagram $(*)$ in Definition \ref{11} is equivalent to 
a DGA homotopy commutative diagram: $(***)$ %
$$\xymatrix{
(\Lambda W,d_Y) \ar@/^2pc/[rr]^{{M(f)}} \ar[r]_{i_f}& ( \Lambda  W\otimes \Lambda U,D)&(\Lambda V,d_X)\ar[l]^{\simeq}\\
(\Q[t_1,..,t_r]\otimes \Lambda W,D_{2}) \ar@/^2pc/@{.>}[rr] \ar@{.>}[r]_{\ \ \ \ \ \   }\ar[u]^{p_t}&(\Q[t_1,..,t_r]\otimes\Lambda  W\otimes \Lambda U,D'_{1})\ar[u]_{}&(\Q[t_1,..,t_r]\otimes \Lambda V,D_{1})\ar[l]^{\ \ \ \ \ \simeq}\ar[u]\\
(\Q[t_1,..,t_r],0)\ar[u]\ar@{=}[r]&(\Q[t_1,..,t_r],0) \ar[u]\ar@{=}[r]&(\Q[t_1,..,t_r],0) \ar[u]
}$$
with $\dim H^*(\Q[t_1,..,t_r]\otimes \Lambda W,D_2)<\infty$ 
and $\dim H^*(\Q[t_1,..,t_r]\otimes \Lambda V,D_1)<\infty$.

For example,  for the fibre inclusion of the Hopf fibration $f:S^3\to S^7$,
$r_0(S^3)=r_0(S^7)=r_0(f)=1$
since it induces the natural inclusion
$E_1=\C P^1\to \C P^3=E_2$
satisfying $(*)$ (without rationalization).
On the other hand, for a rationally non-trivial fibration $S^5\to X\overset{f}{\to} Y=S^3\times S^3$,
$r_0(X)=1$, $r_0(Y)=2$ and $r_0(f)=0$ from $(***)$.
If $r_0(f)=0$,
the map  $f$ can not (rationally) be an $S^1$-equivariant map preserving almost free actions.

 From the definition, 
 $$r_0( f)\leq \min \{  r_0(X),r_0(Y) \}$$
for any map $f:X\to Y$.
In particular, $$r_0(i_X)=r_0(X) \mbox{  \ \ and  \ \ }r_0(p_Y)=r_0(Y)$$
for the inclusion  $i_X:X\to X\times Y$, 
 for the projection $p_Y:X\times Y\to Y$.
 


Recall the LS category ${\rm cat}(f):=\min \sharp 
\{ U_i\subset X\mid X=\cup_iU_i\mbox{ is an open covering with } f|_{U_i}\simeq *\} -1$ for  a map $f:X\to Y$,
where  ${\rm cat}(id_X)={\rm cat}(X)$, the LS category of a space $X$.
Here $\sharp$ denotes the cardinality of a set.
It satisfies ${\rm cat}( f)\leq \min \{  {\rm cat} (X),{\rm cat}(Y) \}$
for any map $f:X\to Y$.

\begin{thm} For maps $f:X\to Y$ and  $g:Y\to Z$,
 $r_0(g\circ f)$
can be arbitrarily large
compared  with $r_0(f)$ and $r_0(g)$.
\end{thm}

This theorem follows from the second example in 
 \begin{ex}\label{kkkk}
{\rm 
(1) For any $m,n$ and $s\leq \min \{m,n\}$,
there are maps
$f:X\to Y$ and $g:Y\to Z$ with $r_0(f)=m$, 
$r_0(g)=n$ and  $r_0(g\circ f)=s$.
For example,
put $$f:S^3_1\times \cdots  \times S^3_m\to S^3_1\times \cdots \times  S^3_m
\times S^5_1\times \cdots  \times S^5_n\ \ \ \mbox{and} $$
$$g: S^3_1\times \cdots  \times S^3_m
\times S^5_1\times \cdots \times S^5_n
\to S^3_1\times \cdots  \times S^3_s\times S^5_1\times \cdots  \times S^5_{n-s} $$
where $f|_{S^3_i}=id_{S^3_i}$ for all $i$,
 $g|_{S^3_i}=id_{S^3_i}$ for $i=1,..,s$,
   $g|_{S^5_i}=id_{S^5_i}$ for $i=1,..,n-s$
and $ g|_{S^n_i}=*$ for other $i$.
Then we have an  example of it.

(2) Consider the maps 
$f:X\to Y$ and $g:Y\to Z=S^3\times \cdots \times S^3$ (2n-factors)  with the following models.
For $n> 1$, 
put $$M(g):M(Z)=(\Lambda (w_1,..,w_{2n}),0)\to (\Lambda (w_1,..,w_{2n},w),d_Y)=M(Y)$$
with $|w_i|=3$ for all $i$,  $|w|=6n-1$, $d_Y(w)=w_1\cdots w_{2n}$
and
$$M(f):M(Y)=(\Lambda (w_1,..,w_{2n},w),d_Y)\to (\Lambda (w_1,..,w_{2n},w,y),d_X)= M(X)$$
with $|y|=5$,  $d_X(w)=w_1\cdots w_{2n}$ and $d_X(y)=w_1 w_2$.
Then we have $$r_0(f)=1,\ r_0(g)=0 \mbox{ and }r_0(g\circ f)=2n-2.$$
In particular  we can verify the third since
$$M(X)=(\Lambda (w_1,..,w_{2n},w,y),d_X)\cong 
(\Lambda (w_3,..,w_{2n},w),0)\otimes A$$
with $A:=(\Lambda (w_1,w_{2},y),d_X)$
induces the $\Q [t_1,..,t_{2n-2}]$-map  
$$F: (\Q [t_1,..,t_{2n-2}]\otimes \Lambda (w_1,..,w_{2n}),D)\to  (\Q [t_1,..,t_{2n-2}]\otimes 
\Lambda (w_3,..,w_{2n},w),D)\otimes A$$
with $Dw_i=t_{i-2}^2$ for $i=3,..,2n$
and $F(w_i)=w_i$ for all $i$.
}
\end{ex}

On the other hand, 
${\rm cat} (g\circ f)\leq \min \{  {\rm cat} (f),{\rm cat}(g) \}$
\cite[Exercise 1.16]{C}.  
 Futhermore,   we know 
 ${\rm cat} (c)=0$ 
for the  constant map $c:X\to Y$ for any space $Y$.
But we can often rationally construct  a suitable model $M(Y)$  
such that $r_0(X)=r_0(c)=r_0(Y)$.
For example, for $M(X)=(\Lambda (x,y,z),d)$
with $|x|=3$,  $|y|=5$,  $|z|=7$,  $dx=dy=0$ and $dz=xy$,
put 
$M(Y)=(\Lambda (x',y',z'),d)$
with $|x'|=5$,  $|y'|=7$,  $|z'|=11$,  $dx'=dy'=0$ and $dz'=x'y'$.
Then we can construct commutative diagram
$$\xymatrix{M(Y)=(\Lambda (x',y',z'),d)\ar[r]^0&(\Lambda (x,y,z),d)=M(X)\\
M(E_2)=(\Q [t]\otimes \Lambda (x',y',z'),D_2)\ar[r]^F\ar[u]^{t=0}&(\Q [t]\otimes \Lambda (x,y,z),D_1)=M(E_1)\ar[u]^{t=0}
}
$$
where $F(t)=t$, $F(x')=xt$,  $F(y')=yt$, $F(z')=zt^2$, $D_1z=xy+t^4$ and  $D_2z'=x'y'+t^6$.
Since  $\dim H^*(E_i;\Q )<\infty$, we see $r_0(X)=r_0(c)=r_0(Y)=1$.
Thus the two numerical invariants of a map, $r_0(f)$ and ${\rm cat}(f)$,
have  very different properties.







 \section{Examples}
 Suppose that  $G$ and $K$ are compact connected  Lie groups and 
 $K$ is a compact connected subgroup of $G$.
Recall the result of Allday-Halperin \cite[Remark(1)]{AH}:

\begin{thm}\label{AH}\rm{(\cite[Corollary 4.3.8]{AP},\cite[Corollaries 7.14 and  7.15]{FOT})} \ 
$$r_0(G)={\rm rank} G{ \ \  and\ \ }r_0(G/K)={\rm rank} G-{\rm rank} K.$$
\end{thm}

Theorem \ref{AH} says that
there is a  pure (two stage) Borel fibration $G/K\to ET^r\times_{T^r}G/K\to BT^r$   
 (\cite{HT}) with $\dim H^*(   ET^r\times_{T^r}G/K ;\Q )<\infty$
for  $r={\rm rank} G-{\rm rank} K$; i.e.,
the differential $D$ in the relative model   of $(***)$ in \S 1 satisfies 
$D_1v\in \Q[t_1,..,t_r]\otimes \Lambda V^{even}$ for $v\in V^{odd}$
and $D_1v=0$ for $v\in V^{even}$
when $M(G/K)=( \Lambda V,d)$.

\begin{thm}\label{P}  Let  $G$ and $K$ be simply connected Lie groups and 
 $K$ a compact connected subgroup of $G$.
For a principal $K$-bundle  $K\overset{g}{\to} G\overset{f}{\to} G/K$, 
  $ r_0(g)={\rm rank} K$
and $ r_0(f)={\rm rank} G-{\rm rank} K$.
\end{thm}

\noindent{\it Proof}.
Put the relative model of $\xi :G\overset{f}{\to} G/K\overset{k}{\to} BK$ as
$$M(BK)=(\Q[x_1,..,x_n],0)\to (\Q[x_1,..,x_n]\otimes \Lambda V,d)\overset{p}{\to} (\Lambda V,0)=M(G)$$
with $dv_i\in \Q[x_1,..,x_n]$ for $v_i\in V$ \cite[Proposition 15.16]{FHT}.
Put $r={\rm rank} G-{\rm rank} K$.
From Theorem \ref{AH} and Proposition \ref{H},
 there is a DGA
$A:=(\Q[t_1,..,t_r]\otimes \Q[x_1,..,x_n]\otimes \Lambda V,D)$
 where $Dv_i=dv_i+g_i$ with $g_i\in (t_1,..,t_r)$, $Dx_i=0$ 
and   $\dim H^*({A})<\infty$.
The DGA-projection $p$ is extended to the $\Q[t_1,..,t_r]$-projection 
$$F:A=(\Q[t_1,..,t_r]\otimes \Q[x_1,..,x_n]\otimes \Lambda V,D)
\to (\Q[t_1,..,t_r]\otimes \Lambda V,\overline{D})=:B,$$
which induces
 $\dim H^*({B})<\infty$. 
Thus $ r_0(f)\geq {\rm rank}  G-{\rm rank} K$.
\hfill\qed\\

\begin{ex}\label{2}{\rm Let $SU(n)$ be the n-th special unitary group.
Then $M(SU(6))$ is given as $(\Lambda (v_1,v_2,v_3,v_4,v_5),0)$
with $|v_1|=3$, $|v_2|=5$, $|v_3|=7$, $|v_4|=9$ and  $|v_5|=11$.

(1) For  the principal bundle $SU(3)\overset{g}{\to} SU(6)\overset{f}{\to} SU(6)/SU(3)$,
the relative model is   extended to
$$\xymatrix{\Lambda (v_3,v_4,v_5),0\ar[dr]_{}\ar[r]\ar@{=}[d]& \Lambda (v_1,v_2,v_3,v_4,v_5),0
\ar[r]^{M(g)} &\Lambda (v_1,v_2),0\\
M( SU(6)/SU(3))& \Q[t_1,t_2]\otimes  \Lambda (v_1,v_2,v_3,v_4,v_5),D_2\ar[r]^{\ \ \  F}\ar[u]&\Q[t_1,t_2]\otimes  \Lambda (v_1,v_2),D_1\ar[u]
}$$
with
 $D_1v_1=D_2v_1=t_1^2$, $D_1v_2=D_2v_2=t_2^3$, $D_2v_3=D_2v_4=D_2v_5=0$.
Thus $ r_0(g)={\rm rank}  SU(3)=2$.

(2) For  the principal bundle $SU(3)\times SU(3)\overset{g}{\to} SU(6)\overset{f}{\to} SU(6)/SU(3)\times SU(3)$, 
the relative model  is extended to
{\small $$\xymatrix{\Q[x_1,x_2,x_1',x_2']\otimes \Lambda V,d\ar[dr]\ar[r]& \Lambda V\otimes C,D\ar[r]& \Lambda (u_1,u_2,u'_1,u'_2),0 \\
M(SU(6)/SU(3)\times SU(3))\ar[u]^{\simeq}& \Q[t_1,..,t_4]\otimes \Lambda V\otimes C,D_t\ar[r]^{F\ \ \ }\ar[u]& \Q[t_1,..,t_4]\otimes\Lambda (u_1,u_2,u'_1,u'_2),\overline{D}_t \ar[u]
}$$}
where $|x_1|=|x_1'|=4$,  $|x_2|=|x_2'|=6$, $V=\Q( v_1,v_2,v_3,v_4,v_5)$, $dx_i=dx_i'=0$,
$dv_1=x_1+x_1'$,
$dv_2=x_2+x_2'$,
$dv_3=x_1^2+{x_1'}^2$,
$dv_4=x_1x_2+x_1'x_2'$, 
$dv_5=x_2^2+{x_2'}^2$
\cite[p.486]{GHV}
and 
$$C=\Q[x_1,x_2,x_1',x_2']\otimes\Lambda (u_1,u_2,u'_1,u'_2) $$
with $|u_i|=|u_i'|=2i+1$, $Du_i=x_i$ and $Du_i'=x_i'$; i.e., $H^*(C)=\Q$. 
Here  $D_tu_1=x_1+t_1^2$,  $D_tu_2=x_2+t_2^3$,  $D_tu'_1=x'_1+t_3^2$,  $D_tu'_2=x'_2+t_4^3$.
Thus $ r_0(g)={\rm rank} SU(3)\times SU(3)=4$.

Also $ r_0(f)={\rm rank} SU(6)-{\rm rank} SU(3)\times SU(3)=1$.
Indeed, for the  minimal model
$M(SU(6)/SU(3)\times SU(3))=(\Q[x_1,x_2]\otimes \Lambda (v_3,v_4,v_5),d)$
 with $dx_1=dx_2=0$, $dv_3=x_1^2$,  $dv_4=x_1x_2$ and   $dv_5=x_2^2$, we have 
a commutative diagram
$$\xymatrix{(\Q[x_1,x_2]\otimes \Lambda (v_3,v_4,v_5),d)\ar[r]^{x_i=0}& (\Lambda (v_1,v_2,v_3,v_4,v_5),0) \\
( \Q[t,x_1,x_2]\otimes \Lambda (v_3,v_4,v_5),D)\ar[r]^{F}\ar[u]^{t=0}& (\Q[t]\otimes\Lambda (v_1,v_2,v_3,v_4,v_5),\overline{D}) \ar[u]_{t=0}
}$$
where    $Dv_3=x_1^2$,  $Dv_4=x_1x_2+t^5$
and   $Dv_5=x_2^2$.
}
\end{ex}

\begin{thm}\label{A}
If  $G/K\overset{g}{\to} X\overset{f}{\to} Y$ is 
a fibration  associated with  a principal $G$-bundle, 
then  $r_0(g)=  {\rm rank} G-{\rm rank} K$.
\end{thm}


\noindent
{\it Proof.} 
Put the model of the fibration  $f:X\to Y$ as $i: M(Y)=( \Lambda W,d_Y)\to 
(\Lambda W\otimes \Lambda V,D)$.
Then $M(G/K)=( \Lambda V,d)$ with $d=\overline{D}$.
Note  that
$Dv\in  \Lambda W\otimes \Lambda V^{even}$ for $v\in V^{odd}$ and   $Dv=0$ for $v\in V^{even}$ \cite[(3.4)]{HT}
from the assumption.
Put $r={\rm rank} G-{\rm rank} K$.
Fix a differential $d_t$ on $\Q[t_1,..,t_r]\otimes  \Lambda V$  with
$\overline{d}_t=d$ and 
  $\dim H^*(\Q[t_1,..,t_r]\otimes  \Lambda V,{d_t})<\infty$,
which exists  from Theorem \ref{AH}.
Note $d_t|_{V^{even}}=0$ and $d_t(V^{odd})\subset \Q [t_1,..,t_r]\otimes \Lambda V^{even}$.
Then we have the  relative model
$$( \Lambda W,d_Y)\to 
(\Q[t_1,..,t_r]\otimes \Lambda W\otimes \Lambda V,D_t)\to 
(\Q[t_1,..,t_r]\otimes  \Lambda V,{d_t})$$
with 
 $D_t(v):=Dv+(d_t-d)(v)$ and $D_t(w):=d_Yw$.
It is embedded into  a  commutative diagram
$$\xymatrix{
( \Lambda W,d_Y)\ar[r]^{  i\ \ }\ar[rd]&(\Lambda W\otimes \Lambda V,D)\ar[r]&
 (\Lambda V,d)=M(G/K)\\
&(\Q[t_1,..,t_r]\otimes \Lambda W\otimes \Lambda V,D_t)\ar[u]_{t_i=0}\ar@{.>}[r]^{\ \ F}&
(\Q[t_1,..,t_r]\otimes \Lambda V,{d_t}).\ar[u]_{t_i=0}
}$$
Since $\dim H^*(\Q[t_1,..,t_r]\otimes  \Lambda V,d_t)<\infty$, we have 
$\dim H^*(\Q[t_1,..,t_r]\otimes \Lambda W\otimes \Lambda V,D_t)<\infty$
from the Serre spectral sequence.
Thus we have $r_0(g)\geq r$.
From Theorem \ref{AH}, we have $r_0(g)= r$.
\hfill\qed\\

\begin{rem}{
If  a fibration $C\overset{g}{\to} X\overset{f}{\to} Y$  is not (associated with) a principal bundle, 
it does not hold that $r_0(g)=r_0(C)$.
For example, for 
the  rational fibration
$SU(3)\to X\to S^3$
given by 
$$(\Lambda w,0)\to (\Lambda (w,u_1,u_2),D)\to (\Lambda (u_1,u_2),0)$$
where $|w|=|u_1|=3$ and $|u_2|=5$ and 
$Du_2=wu_1$, we have 
$r_0(g)=1<2=r_0( SU(3))$.
Also for a fibration over $S^3$  
of the relative model
$$(\Lambda w,0)\to (\Lambda (w,v_2,v_3, v_4,v_5,v_7),D)\to (\Lambda (v_2,v_3, v_4,v_5,v_7),\overline{D})=M(C)$$
where $|w|=3$, $|v_i|=i$, $Dv_2=0$,
$Dv_3=v_2^2$, $Dv_4=wv_2$, $Dv_5=v_2v_4-wv_3$ and $Dv_7=v_4^2+2wv_5$,
 we can check
$r_0(g)=r_0(X)=0$ from \cite{Y1}.
On the other hand, $r_0(C)=1$
since $\dim H^*( \Q [t]\otimes \Lambda (v_2,v_3, v_4,v_5,v_7),\overline{D}_t)<\infty$
by  $\overline{D}_t(v_3)=v_2^2$, $\overline{D}_t(v_5)=v_2v_4+t^3$ and $\overline{D}_t(v_7)=v_4^2$. 
Note  their  fibres are  pure but the above  two fibrations are  not pure \cite{HT}.  
Compare with Theorem \ref{A}.
}\end{rem}


\begin{thm}\label{em}
For an odd-spherical   fibration  $\xi :S^{2n-1} \to X\overset{f}{\to} Y$, 
 suppose 
$\pi_{>2n}(Y)\otimes \Q=0$. 
Then,  for  any  free $T^r$-action $\mu$ on $X'$ with $X'_{\Q}\simeq X_{\Q}$
such that  $r_0(ET^r\times_{T^r}^{\mu}X')=0$,
 there is no  map $F$
between fibrations 
$$\xymatrix{
X_{\Q} \ar[r]\ar[d]_{f_{\Q}}&(ET^r\times_{T^r}^{\mu}X')_{\Q}\ar[r]\ar[d]^{F}&BT^r_{\Q}\ar@{=}[d]\\
Y_{\Q}\ar[r]&(ET^r\times_{T^r}^{\tau}Y')_{\Q}\ar[r]&BT^r_{\Q}\\
}$$
such that 
 $\tau$ is  a free $T^r$-actionon $Y'$ with $Y'_{\Q}\simeq Y_{\Q}$.
 In particular, $r_0(f)<r_0(X)$.
\end{thm} 

 \noindent
 {\it Proof.}
Put $M(S^{2n-1})= (\Lambda y,0)$ and $M(Y)=(\Lambda V,d_Y)$.
Supppose that there is a map     $(\Q [t_1,..,t_j]\otimes \Lambda V,D_2)
\to (\Q [t_1,..,t_j]\otimes \Lambda V\otimes \Lambda y,D_1)$
with $\dim H^*(D_i)<\infty$.
Then, from degree reasons,  there is a KS-extension of  $(\Q [t_1,..,t_j]\otimes \Lambda V\otimes \Lambda y,D_1)$  by $t_{j+1}$, 
 the DGA   $A:=(\Q [t_1,..,t_{j+1}]\otimes \Lambda V\otimes \Lambda y,D')$ with $$D'(y)=D_1(y)+t^{n}_{j+1}
\mbox{\ \  \ and\ \ \   }D'(v)=D_1(v)\ \mbox { for }v\in V$$  
satisies  $\dim H^*(A)<\infty$.
  \hfill\qed\\
  
\begin{qu}\label{yamaguchi} For a fibration $\xi :C\overset{g}{\to} X\overset{f}{\to} Y$ 
with  fibre $C$   simply connected   
of  finite rational cohomology, does it hold that 
$r_0(g)+r_0(f)\leq r_0(X)$ ?
\end{qu}



\begin{rem}\label{product}{
The above question  is true for many cases.
For example, it is true 
for  
the fibrations of Theorem \ref{em}
or when $r_0(C)=0$.
Of course, 
it is true when $\xi$ is rationally trivial.
But it may not be equal.  
Recall 
  Halperin's  inequality $r_0(X)=r_0(X)+r_0(S^{2n})<r_0(X\times S^{2n})$   for a
formal  space $X$ and an integer $n>1$ \cite{JL}. 
 For any even integer $n$, there is a space $X_n$ such that
$r_0(X_n)=0$ and 
$r_0(X_n\times S^{6n+1})\geq n$.
In the following, we give an example of the  model.
Put 
 $$M(X_n)=(\Lambda V,d)=(\Lambda (u_1,u_2,..,u_n,v_1,v_2,...,v_{2n},v,w),d)$$ with
$$dv_i=du_i=dw=0\ \  \mbox{ for all $i$ and}$$
$$dv=
u_1u_2u_3\cdots u_n(v_1v_2+v_3v_4+v_5v_6+\cdots +v_{2n-1}v_{2n})+w^2$$
and $|v_i| =|u_i|=3$ for all $i$, $|w|=(3n+6)/2$,   $|v|=3n+5$.
Then we can check that $r_0(X_n)=0$
by Proposition \ref{H}
since $\dim H^*(\Q [t]\otimes \Lambda V,D)=\infty$
for any differential $D$ by direct calculations.
Put $M(S^{6n+1})=(\Lambda y,0)$ with $|y|=6n+1$ and
$$M(ET^n\times _{T^n}(X_n\times S^{6n+1}))=(\Q [t_1,..,t_n]\otimes \Lambda V\otimes \Lambda y ,D)$$ by
$$Dv=dv
+\sum_{i=1}^nu_iyt_i,\ \ Dv_{2i}=t_i^2, \ \ Dv_{2i-1}=0\ (i=1,..,n)$$
$$\mbox{and }\ \ Dy=\sum_{i=1}^n(-1)^{i+1}v_{2i-1}u_1\cdots \hat{u_i}\cdots u_{n}t_i.$$
Then $D\circ D=0$ and $\dim H^*(\Q [t_1,..,t_n]\otimes \Lambda V\otimes \Lambda y )<\infty$.
Thus  $r_0(X\times Y)$ can be arbitrarily large compared to  
 $r_0(X)+r_0(Y)$.
}
\end{rem}

\begin{rem}{  
Is there a good cohomological upper bound for $r_0(f)\ ?$ 
Recall that
S.Halperin proposes  the {\it toral rank  conjecture} (TRC) that
the inequality $$\dim H^*(X;\Q)\geq 2^{r_0(X)}$$ holds \cite{H} (\cite[{\bf 39}]{FHT},
\cite[Conjecture 7.20]{FOT}).
For example, a homogeneous space satisfies it \cite[(7.23)]{FOT}.
It is natural to ask whether
the inequality
$\dim Im (H^*(f;\Q))\geq 2^{r_0(f)}$
holds.
But that is not  the case in general.
For example, put 
$M(X)=(\Lambda (v_1,v_2,v_3),0)$ with $|v_i|=3$
and 
$M(Y)=(\Lambda (x,y,v_1,v_2,v_3),d)$
with $dv_1=x^2$, $dv_2=xy$, $dv_3=y^2$, $dx=dy=0$, $|x|=|y|=2$,
and $M(f)(v_i)=v_i$ and $M(f)(x)=M(f)(y)=0$.
Then $H^*(f;\Q)$ is trivial; i.e., $ \dim Im (H^*(f;\Q))=1$. 
On the other hand,  $r_0(f)=1$.
Indeed,
$(\Q[t]\otimes \Lambda (x,y,v_1,v_2,v_3) ,D_{2})$
is given by $D_{2}v_1=x^2$,
$D_{2}v_2=xy+t^2$,
$D_{2}v_3=y^2$
and 
$(\Q[t]\otimes \Lambda (v_1,v_2,v_3),D_{1})$
is given by
 $D_{1}v_1=D_{1}v_3=0$,
$D_{1}v_2=t^2$.
Then 
$\dim H^*(D_{1})<\infty$, $\dim H^*(D_{2})<\infty$
and 
$M(f)$ is extended to a $\Q[t]$-morphism $F: (\Q[t]\otimes \Lambda (x,y,v_1,v_2,v_3) ,D_{2})\to (\Q[t]\otimes \Lambda (v_1,v_2,v_3),D_{1})$
with $F(x)=F(y)=0$.



}
\end{rem}





\end{document}